\documentclass{amsart}

\usepackage{fullpage}
\usepackage{amssymb}
\usepackage{amsmath}
\usepackage{amsfonts}
\usepackage{amsxtra}
\usepackage{amscd}
\usepackage{epsfig}
\usepackage{enumerate}
\usepackage{float}
\usepackage{xy}

\xyoption{all}

\theoremstyle{plain} \newtheorem{Thm}{Theorem}[section]
\theoremstyle{plain} \newtheorem{Cor}[Thm]{Corollary}
\theoremstyle{plain} 
\theoremstyle{plain} 
\theoremstyle{definition} 
\theoremstyle{definition} \newtheorem{Exa}[Thm]{Example}
\theoremstyle{remark} \newtheorem{Rem}[Thm]{Remark}
\newenvironment{Proof}{{\bf Proof :}}{$\Box$\\}

\newcommand{\newfontobj}[2]{
  \newcommand{#1}[1]{
    \expandafter\def\csname##1\endcsname{{#2 ##1}}}}

\newfontobj{\class}{\rm}
\newfontobj{\lang}{\bf}
\newfontobj{\oper}{\rm}

\class{PSPACE}		
\class{AM}		
\class{MA}
\class{NP}
\class{UP}
\class{P}
\class{RP}
\class{BPP}
\class{DTIME}
\class{ZPP}
\class{EXPSPACE}
\class{coNP}
\class{coRP}
\class{coAM}
\class{PH}
\class{co}
\class{BP}	

\lang{HN} 		
\lang{SAT}
\lang{PIM}
\lang{CSG}
\lang{ESC}
\lang{IP}
\lang{PRIMES}

\oper{Pr}    
\oper{E}     
\oper{Ord}   
\oper{li}

\floatstyle{ruled}
\newfloat{Algorithm}{thp}{loa}[section]

\title{The Hardness of Computing an Eigenform}
\author{Eric Bach}
\address{Department of Computer Science, University of Wisconsin-Madison, Madison, WI - 53706.}
\email{bach@cs.wisc.edu}

\author{Denis Charles}
\address{Microsoft Research, One Microsoft Way, Redmond, WA - 98052}
\email{cdx@microsoft.com}

\date{8 September, 2006}

\begin{document}
\maketitle
\section{Introduction}
The Fourier coefficients of modular forms encode very interesting arithmetic data. For example,
divisor sums, partition numbers, trace of Frobenius of the reduction modulo primes of an elliptic curve over $\mathbb{Q}$, and more generally,
trace of Frobenius of many Galois representations of dimension $2$ over finite fields (this being a conjecture of
Serre) are all known to be, or conjectured to be, Fourier coefficients of modular forms. A particularly important
family of modular forms are the so-called Hecke eigenforms. 
These are modular forms that are also simultaneous eigenforms for
an algebra of operators called the Hecke operators that operate on the spaces of modular forms. The Fourier
coefficients of Hecke eigenforms are particularly important arithmetically. 
Indeed, many of the examples given above arise as Fourier coefficients
of Hecke eigenforms. \\

In this article we are concerned with the computational
complexity of computing the Fourier coefficients of these Hecke eigenforms. Currently, there are three approaches
to computing the Fourier coefficients of modular forms: a method based on computing theta series of lattices
\cite{piz80}; the method of modular symbols \cite{mer94, ste00}; and one 
based on the Selberg Trace formula (see \cite{comp_tau} and Chapter 5 of \cite{my_thesis}). All of these methods result in algorithms with
exponential running time to compute the Fourier coefficients. The Fourier coefficients of Hecke eigenforms are
multiplicative and satisfy recurrences for prime powers. Since there
are subexponential time algorithms for factoring integers, the interesting problem is to compute
the $p$-th Fourier coefficient, for prime $p$, efficiently. However, this problem is still open in general. For any fixed
eigenform of weight $2$ one can use Schoof's algorithm (\cite{sch85}) for counting points on elliptic curves over finite fields
to compute the $p$-th Fourier coefficient efficiently. Recent work of Edixhoven et al. suggests that this approach generalizes to compute eigenforms of weight $k \geq 2$ (\cite{edix06}). \\

There have been no hardness results known for computing the
Fourier coefficients of eigenforms (except for those of Eisenstein series where the hardness results follow from \cite{bms86}). 
In this article, we give evidence that computing Fourier coefficients
of the Hecke eigenforms for composite indices is no easier than factoring integers. More precisely, we show that
the existence of a polynomial time algorithm that, given $n$, computes the $n$-th Fourier coefficient of a (fixed) Hecke eigenform
implies that we can factor most RSA moduli (numbers that are products of two distinct primes) in polynomial time. In particular, our
result implies a hardness result for computing the Ramanujan tau function.

\subsection{Preliminaries and notation}
Since there are a number of excellent references for modular forms we refrain from reproducing the definitions here.
Instead, we refer the reader to any of the books 
\cite{ser70, shi71, lan76, kob93, ono04, ds05} for the definition and background on modular forms. 
The Fourier expansions
of modular forms that we refer to are the Fourier expansions at the cusp $\infty$. In the Fourier expansions $q$ stands
for $e^{2\pi \imath z}$. The letters $p$ and $q$ will be used for (rational) prime numbers,
the latter will be used when no confusion can arise with the $q$-expansions.

\section{The Reduction}

Let $S_k(\Gamma_0(N),\chi)$ be the space of cusp forms of {\em even\footnote{The condition of even weight is needed for a technical reason in the proof.} integer} weight $k$ ($\geq 2$), level $N$, and character $\chi ~\mathrm{mod}~ N$. 
In what follows, fix $f(z) = \sum_{1 \leq n} a(n) q^n \in S_k(\Gamma_0(N),\chi)$ 
to be a
normalized ($a(1) = 1$) Hecke eigenform. We will also assume that $f$ is not of CM type in the
sense of Ribet \cite{rib77}. This means that there does not exist an imaginary quadratic field, $L$, such
that $a(p) = 0$ for all primes $p$ that are inert in $L$. Under
these assumptions, a beautiful theorem of Serre (\cite{ser81} Corollary 2 to Theorem 15) gives
us bounds on the number of primes $p$ for which $a(p) = 0$.

\begin{Thm}\label{thm_density_zero:this}
Let $f(z) = \sum_{1 \leq n} a(n) q^n \in S_k(\Gamma_0(N),\chi)~ (k \geq 2)$ be
a normalized Hecke eigenform that is not of CM type.
Define $P_{f}(x) = \sharp\{p \leq x~:~ p \text{ a prime such that } a(p) = 0\}$. Then
\begin{align*}
P_f(x) = O\left( \frac{x}{(\log x)^{\frac{3}{2}-\delta}}\right) \text{ for all } \delta > 0.
\end{align*}
Moreover, if one assumes the Generalized Riemann Hypothesis, we have
\begin{align*}
P_f(x) = O\left( x^{\frac{3}{4}}\right).
\end{align*}
\end{Thm}
The assumption that $f$ not be of CM type is necessary, see Remark \ref{rem_cmtype:this}.\\
 
The Fourier coefficients of a normalized Hecke eigenform need not be integers, but they are at least
 algebraic integers (see \cite{ono04} \S2.4 \& \S2.5; the result also follows from \cite{shi71} Theorem 3.52).
Furthermore, we know that each eigenvalue lies in a number field of degree at most $\dim S_k(\Gamma_0(N),\chi) \varphi(N)$
since the characteristic polynomials of the Hecke operators have degree $\dim S_k(\Gamma_0(N),\chi)$
over the field $\mathbb{Q}(\chi)$.
In fact, the field $K_f = \mathbb{Q}(a(2),a(3),\cdots,a(n),\cdots)$ is a number field and so a finite degree extension
of $\mathbb{Q}$. Since $f$ is fixed we can assume that we can do computations in this field efficiently. 
We assume that the supposed algorithm that computes the Fourier coefficients
takes as input an integer $n$ and gives us the (monic) minimal polynomial of the $n$-th Fourier coefficient $a(n)$. We also assume that the algorithm provides a complex approximation to $a(n)$ that distinguishes $a(n)$ from its conjugates. In other words, we not only have the sub-field $\mathbb{Q}(a(n)) \subseteq K_f$,
but we also have an embedding of $\mathbb{Q}(a(n))$ in $\mathbb{C}$.
Since the space $S_k(\Gamma_0(N),\chi)$ is fixed, and $\chi$ is a Dirichlet
character $\mathrm{mod}~ N$, we can also compute $\chi(n)$ for any integer $n$.\\

Next, we describe how we can factor RSA moduli if we can compute the Fourier coefficients $a(n)$. 
We are given a positive integer $n = pq$, where $p, q$ are distinct odd primes. We
can also assume (without loss of generality) that $\gcd(N,n) = 1$. 
Let $x = a(p)/\chi(p)p^{\frac{k-1}{2}}$ and 
$y = a(q)/\chi(q)q^{\frac{k-1}{2}}$. Note that we can also assume that $\chi(n) \neq 0$
for otherwise $\gcd(n,N) \neq 1$. We will make the assumption that $x\neq 0$ and $y \neq 0$ in the following analysis.\\

Using the algorithm to compute the Fourier coefficients of $f$ we can compute
\begin{align*}
	A &=_{\rm def} n^{\frac{k-1}{2}} xy \\
	  &= \frac{a(n)}{\chi(n)}
\end{align*}
and
\begin{align*}
	B =_{\rm def} a(n^2).
\end{align*}

Now by multiplicativity and the recurrences for prime powers that $a(n)$ satisfy we have (\cite{kob93} III.\S5)
\begin{align}
B = a(n^2) &= a(p^2)a(q^2) \\
	   &= (a(p)^2 - p^{k-1}\chi(p))(a(q)^2 - q^{k-1}\chi(q)) \\
	   &= n^{k-1} \chi(n)(x^2 - 1)(y^2 - 1). \label{eqn_square_coeff:this}
\end{align}

Thus we have a pair of simultaneous equations for $x$ and $y$ which we can solve.
Setting $\alpha = A/n^\frac{k-1}{2}$ and $\beta = B/n^{k-1} \chi(n)$, one obtains
\begin{align*}
x^2 &= \frac{\alpha^2 - \beta + 1 \pm \sqrt{(\alpha^2 -\beta + 1)^2 - 4\alpha^2}}{2}
\end{align*}
and
\begin{align*}
y^2 &= \frac{\alpha^2}{x^2}.
\end{align*}
Substituting the definitions of $\alpha$ and $\beta$ and clearing denominators we get
\begin{align*}
x^2 = \frac{\left(A^2\chi(n) - B + n^{k-1}\chi(n)\right) \pm 
\sqrt{\left( A^2\chi(n) - B + n^{k-1}\chi(n)\right)^2 - 4 A^2 n^{k-1}\chi(n)^2}}{2 \chi(n) n^{k-1}}.
\end{align*}

We note that the radicand is the square of an algebraic integer (see below) and hence the square root can be
computed exactly. This can be computed efficiently by computing numerical approximations
to the square roots of all the conjugates of radicand. By the definition of $x$ we have that
\begin{align*}
x^2 = \frac{a(p)^2}{\chi(p)^2 p^{k-1}}.
\end{align*}
Note that this quantity is not zero under our assumption that $x \neq 0$. We claim that $x^2$ cannot be an algebraic integer if $p$ is large enough. For otherwise, since $k-1$ is odd, this
would make $\sqrt{p}$ an element of $\mathbb{Q}(\chi,a(2),a(3),\cdots,a(n),\cdots)$, but the latter is a finite extension
and thus if $p$ is large enough it cannot contain $\sqrt{p}$. Thus we can recover $p$ from the above expression by
taking the $\gcd$ of the denominator of the above expression with $n$. Since the quantity is
an algebraic number the (reduced) denominator in the expression is the leading coefficient of the minimal polynomial over $\mathbb{Z}$.\\

Suppose $x = 0$ but $y \neq 0$ (i.e. $a(p) = 0$ but $a(q) \neq 0$), we can still proceed as follows. By equation 
(\ref{eqn_square_coeff:this}) we find that $B = n^{k-1} \chi(n)(1 - y^2)$. Thus we can still get $y^2$
and by the above argument find $q$. \\

Thus our reduction will succeed in factoring the integer $n$, unless both $a(p)$ and $a(q)$ are zero.
Since the set of such primes is density $0$ (by Theorem \ref{thm_density_zero:this}), we get the following theorem:

\begin{Thm}\label{thm_eigenhard:this} Let $f(z) = \sum_{1 \leq n}a(n) q^n\in S_{2k}(\Gamma_0(N),\chi)$ 
be a normalized Hecke eigenform that is not of CM-type. Suppose there is a polynomial time algorithm that
computes $a(n)$ given $n$. Then there is a polynomial time algorithm that factors a density $1$ subset
of the RSA moduli.
\end{Thm}

In the case that $f \in S_{k}(\Gamma_0(N),\chi)$ and $k$ is {\em odd} the entire reduction works as long as $p^{\frac{k-1}{2}}$
does not divide $a(p)$ for one of the primes dividing $n$. The failure of the reduction occurs very rarely. 
Indeed, if $k \geq 3$
and $k$ is odd then this implies that $a(p) \equiv 0 \mod p$ which means that $p$ 
is a, so called, {\em non-ordinary prime.} A heuristic
argument given in \cite{gou97} shows that the number of non-ordinary primes below $x$ is $O(\log \log x)$.
Thus, it is likely that the result of Theorem \ref{thm_eigenhard:this} remains true even for odd weight
cuspidal eigenforms.

\begin{Exa} We illustrate the reduction in the case of the Ramanujan Tau function $\tau(n)$ that gives
the Fourier coefficients of $\Delta$, a weight $12$ eigenform of level $1$ and trivial character 
(see \S\ref{sec_tau:this}).\medskip

Let $n = 15$, from the tables in \cite{leh43} one sees that $\tau(15) = 1217160$ and $\tau(15^2) = 2897808426675$.
In the notation of the proof of the theorem we have
\begin{align*}
\alpha^2 &= \frac{81288256}{474609375}, \text{ and } \\
\beta &= \frac{1431016507}{4271484375}.
\end{align*}
From this one finds that $x^2 = \frac{933156}{1953125}$, and $\gcd(1953125,15) = 5$.
\end{Exa}

\begin{Exa} The space $S_{4}(\Gamma_0(29))$ with trivial character has a newform, $f$ (say),
whose expansion begins $q + \gamma q^2 + (-3\gamma -8) q^3 + (-2\gamma -7)q^4 + (4\gamma-1)q^5 + \cdots$, where
$\gamma$ is a root of $x^2 + 2x -1$. A short computation in MAGMA (\cite{magma}) tells us that the $15$th Fourier coefficient is
$-5\gamma-4$ and that the $225$th Fourier coefficient is $-2680\gamma - 6168$. MAGMA computes
that either
\begin{align*}
x^2 = \frac{1}{27}(30\gamma + 73) \text{ or }
x^2 = \frac{1}{125}(-40\gamma+17)
\end{align*}
corresponding to the two square roots of $(\alpha^2-\beta+1)^2-4\alpha$ (again we have
preserved the notation used in the proof). In any case, the denominators
in these expressions yield a proper factor of $15$.
\end{Exa}

\begin{Rem} The hardest cases of factoring the RSA moduli are believed to be those 
of the form $pq$ where the primes $p$ and $q$
are both approximately the same size. One might wonder if the set of RSA moduli on which our reduction works
includes such numbers also. This is indeed true. The number of RSA moduli below a bound $x$ that have both
the factors being approximately the same size can be estimated as follows. Let $c$ be a constant
with $0 < c < 1$. Using the prime number theorem
the number of RSA moduli $pq$ for which $c \sqrt{x} \leq p,q \leq \sqrt{x}$ is 
\begin{align*}
\binom{\pi(\sqrt{x}) - \pi(c\sqrt{x})}{2} = \Theta\left( \frac{x}{\log^2 x}\right), 
\end{align*}
where $\pi$ is the prime counting function. Meanwhile, the number of RSA moduli with primes of the same size
for which our reduction fails is bounded above by (using Theorem \ref{thm_density_zero:this})
\begin{align*}
\binom{\frac{\sqrt{x}}{\log^{3/2 - \delta}x}}{2} = O\left( \frac{x}{\log^{3-2\delta} x}\right), \text{ for all } \delta > 0,
\end{align*}
and if we assume the GRH this upper bound can be strengthened to $O(x^{3/4})$. Thus, our reduction does indeed
work on a density $1$ subset of the ``interesting'' RSA moduli.
\end{Rem}

\begin{Rem} \label{rem_cmtype:this} 
For CM-forms, the prime indexed Fourier coefficients vanish for, roughly, half the primes. And our reduction
will fail if both the prime factors of $n$ are divisible by such primes. 
\end{Rem}

\subsection{Computing a basis of cusp forms}
Theorem \ref{thm_eigenhard:this} has the following consequence
for the problem of computing {\em any} basis of cusp 
forms (with algebraic Fourier coefficients) for $S_k(\Gamma_0(N),\chi)$.

\begin{Cor}\label{cor_basis:this}
Fix $N$, a positive integer, $k \geq 2$, an even integer, and $\chi$ a Dirichlet
character modulo $N$. Assume that $S_k(\Gamma_0(N),\chi)$ contains a Hecke eigenform of non-CM type. Fix also
a basis given by the Fourier expansion
\begin{align*}
	f_i &= \sum_{1 \leq m} a_{i}(m) q^{m}\text{ for } 1 \leq i \leq d,
\end{align*}
such that the $a_{i}(m)$'s are algebraic. Suppose there is a polynomial time algorithm that,
given $n$, computes the list of Fourier coefficients, $a_{i}(m)$, for $1 \leq i \leq d$, then
there is a polynomial time algorithm that can factor a density $1$ subset of the RSA moduli.
\end{Cor}
\begin{Proof}
By our assumption there is a Hecke eigenform not of CM-type in $S_k(\Gamma_0(N),\chi)$. This form
can be normalized by taking a scalar multiple, call this normalized eigenform $g$. 
 Now, since $g$ belongs
to $S_k(\Gamma_0(N),\chi)$ and the $f_i$ span the space we must have that
 $g = \sum_{1 \leq i \leq d} \alpha_i f_i$, where $\alpha_i$ are algebraic numbers.
 The $n$-th Fourier coefficient of $g$ is $\sum_{1 \leq i \leq d} \alpha_i a_{i}(m)$, and
 so this can be computed (in polynomial time) using the supposed algorithm for computing the $a_{i}(m)$'s.
 The result now follows from Theorem \ref{thm_eigenhard:this}.
\end{Proof}

We now investigate the conditions under which the assumption made
in Corollary \ref{cor_basis:this} (that $S_k(\Gamma_0(N),\chi)$ contains an eigenform
of non-CM type) holds. A construction due to Hecke \cite{hec37} (also described by Shimura) shows how one can
obtain essentially all the eigenforms of CM-type (see \cite{rib77} \S3). This construction
together with dimension formulas for $S_k(\Gamma_0(N),\chi)$ can be used to show the
existence of eigenforms of non-CM type. 
The results of Theorem 3.5 and Corollary 3.5 of \cite{rib77} summarize the
construction of CM forms by Hecke. Essentially, these results imply that one gets
CM-forms corresponding to quadratic imaginary fields of discriminant $D$ dividing the level $N$, and
each ideal class character of the orders of discriminant $N$ in these fields. From this observation
and bounds on class numbers of imaginary quadratic fields, we find that the number of eigenforms
of CM-type in $S_k(\Gamma_0(N),\chi)$ is bounded above by $N^{\frac{1}{2}+\epsilon}$ for every $\epsilon > 0$.
The dimension of $S_k(\Gamma_0(N),\chi)$ (see \cite{coe77}) on the other hand is $\Omega(kN)$. Furthermore,
the space $S_k(\Gamma_0(N),\chi)$ has a basis of eigenforms; thus, if $N$ and $k$ are large enough there
will always be eigenforms in $S_k(\Gamma_0(N),\chi)$ that are not of CM-type. In other words, for large enough
$k$ and $N$, the assumption made in Corollary \ref{cor_basis:this} holds. Consequently, computing a basis
for such spaces is at least as hard as factoring RSA moduli.

\subsection{The Ramanujan Tau function} \label{sec_tau:this}
The Ramanujan Tau function $\tau(n)$ is defined to be the $n$-th Fourier coefficient
of the Discriminant function $\Delta(z)$ :
\begin{align*}
\Delta(z) &= q \prod_{1 \leq n}(1 - q^n)^{24}, \\
	&= 1 - 24 q^2 + 252 q^3 - 1472 q^4 + 4830 q^5 - 6048 q^6 - 16744 q^7 + \cdots\\
	&=_{\rm def} \sum_{1 \leq n} \tau(n) q^n.
\end{align*}
It is a fact that $\Delta$ is a Hecke eigenform of weight $12$ and level $1$. There are no CM forms of level 1 (since the discriminant
of the underlying CM field must divide the level), 
so $\Delta$ is not a CM form. Moreover, a conjecture of Lehmer states that $\tau(n)$ is never
zero. If we assume Lehmer's conjecture then the proof of our result now yields a slightly stronger conclusion:

\begin{Cor} Assuming Lehmer's conjecture, computing the Ramanujan tau function is at least as hard as factoring
RSA moduli.
\end{Cor}

{\bf Acknowledgement: } The authors would like to thank Tonghai Yang for pointing them to Hecke's construction
of CM forms. 

\newcommand{\etalchar}[1]{$^{#1}$}

\end{document}